\documentclass[12 pt]{report}
\usepackage{amssymb,latexsym,epsfig,mathrsfs,graphpap,graphics,amsmath,amssymb}
\usepackage{amsmath}

\begin{document}

\vspace{.2in}\parindent=0mm

\begin{flushleft}
{\bf\Large { Biquaternion Windowed Linear Canonical  Transform }}

\parindent=0mm \vspace{.2in}
{\bf{  Owais Ahmad$^{1\star},$ Aijaz Ahmad Dar$^{2}$}}
\end{flushleft}

\parindent=0mm \vspace{.1in}
{{\it $^{1\star}$Department of  Mathematics,  National Institute of Technology, Hazratbal, Srinagar -190 006, Jammu and Kashmir, India. E-mail: $\text{siawoahmad@gmail.com}$}}
\footnote{$\star$ Corresponding Author}

\parindent=0mm \vspace{.1in}
{{\it $^{2}$Department of  Mathematics,  National Institute of Technology, Hazratbal, Srinagar -190 006, Jammu and Kashmir, India.E-mail: $\text{daraijaz9999@gmail.com}$}}

\parindent=0mm \vspace{.2in}
{\bf{Abstract:}} In this paper, we introduce the notion of windowed linear canonical transform in biquaternion setting namely Biquaternion Windowed Linear Canonical Transform (BiQWLCT) and  various properties of BiQWLCT, such as linearity, shift, parity, orthogonality relation, inversion formula, Plancherel theorem are established.  Heisenberg uncertainty principle  associated with the 
Biquaternion Windowed Linear Canonical Transform is also derived. Towards the culmination, an example and some potential applications are presented.

\parindent=0mm \vspace{.1in}

\parindent=0mm \vspace{.1in}

\parindent=0mm\vspace{0.4in}
{\bf{1.  Introduction }}

\parindent=0mm \vspace{.in}
One of the well-known applications of the classical windowed Fourier transform (WFT) is its simultaneous provision of information in the frequency and time domains. It is sometimes referred to as the Gabor transform or the short-time Fourier transform (FT) if a Gaussian window is used. The WFT has been thoroughly studied in several fields, such as quantum physics and communication theory \cite{pp,pq}. The windowed fractional Fourier transform (WFRFT), a generalization of the WFT in the fractional Fourier domain, has been studied in \cite{28,zz}.
The extension of  WFT to the windowed linear canonical transform (WLCT)  has been studied in \cite{24}. This generalized transform is obtained by substituting the LCT kernel for the FT kernel in the WFT formulation.

\parindent=8mm \vspace{.1in}
The linear canonical transform (LCT) is a three-parameter class of linear integral transforms that includes well-known unitary transformations, including the Fourier, fractional Fourier, and Fresnel transforms \cite{13,1}. Because of its additional degrees of freedom and simple geometrical manifestation, the LCT is more flexible than other transforms, making it suitable for investigating complex problems in signal processing, quantum physics, and optics, in general \cite{15,2}. 
However, because of its global kernel, the LCT is unable to display the local LCT frequency contents. The generalized windowed function has been effectively studied using the LCT in  \cite{24,25}. The windowed linear canonical transform (WLCT) is designed to analyze signals with time-varying spectral content. Some significant WLCT features are described, as seen in Refs. \cite{26,24}. These include sampling formulae, the Paley-Weener theorem, uncertain relations, and the equivalent of the Poisson summation formula. Gao and Li \cite{14} introduced the notion of windowed linear canonical transform in quaternion settings and established its various properties and uncertainty relations.

\parindent=8mm \vspace{.1in}
Biquaternions are sometimes known as quaternions with complex components or complex numbers with quaternion real and imaginary parts \cite{3}. This is yet another hypercomplex algebra with bivector, pseudoscalar, vector, and scalar components \cite{4}. Hamilton \cite{3} discovered biquaternions, a generalization of quaternions, in 1853. Ward \cite{5} introduced biquaternion properties. Bekar and Yayh \cite{6} introduced biquaternions (complexified quaternions), split quaternions, and involutions/anti-involutions. This study \cite{7} provided the essential aspects of biquaternions, including a variety of representations.  Sangwine \cite{8} found the roots of $-1$ in the biquaternions. A biquaternion-valued signal, a quaternion with complex components, contains scalar, pseudoscalar, vector, and bivector components \cite{9}. Biquaternion-valued signals have these properties, which allow them to process simultaneous signals with different geometric qualities \cite{4}. Based on their investigation into the extension of the Fourier transform of biquaternion-valued functions, Said et al.\cite{4} introduced the Biquaternion Fourier Transform (BiQFT). They investigated the  relationship between the geometric features of a biquaternion-valued function and the Hermitian symmetries of the biquaternion-valued function (BiQFT).
 The biquaternion Z transformation, proposed lately by Bi et al.\cite{10}, has been used to solve a class of biquaternion recurrence relations. As an extension of the quaternion Fourier transform (QFT), the BiQFT has attracted increasing attention \cite{11,12,14}.
 
 \parindent=8mm\vspace{0.1in}
 Recently, Gao and Li \cite{20} introduced linear canonical transform in biquaternion settings.
 To date, an  attempt has yet to be  made to extend the concept of windowed linear canonical transform in biquaternion settings. Motivated and inspired by the work of Gao and Li \cite{14,20}, we in this paper introduce the notion of windowed linear canonical transform in biquaternion setting namely Biquaternion Windowed Linear Canonical Transform (BiQWLCT) and  various properties of BiQWLCT, such as linearity, shift, parity, orthogonality relation, inversion formula, Plancherel theorem are established.  Heisenberg uncertainty principle  associated with the Biquaternion Windowed Linear Canonical Transform is also derived. Towards the culmination, an example and some potential applications are presented.

\parindent=8mm\vspace{0.1in}
The paper is structured as follows:  In section 2, we recall some  definitions and fundamental characteristics of quaternions and biquaternions. In Section 3, we introduce  the definition of Biquaternion Windowed Linear Canonical Transform (BiQWLCT)  and establish  various properties such as linearity, shift, parity, orthogonality relation, inversion formula and Plancherel theorem.  Heisenberg's uncertainty principle  associated with the 
Biquaternion Windowed Linear Canonical Transform is derived in section 4.  An example and some potential applications are presented respectively in section 5 and section 6.

\parindent=0mm\vspace{0.2in}
{\bf{2. Preliminaries}}

\parindent=0mm\vspace{0.in}
We examine some fundamental quaternion and biquaternion information in this section. Quaternions with complex coefficients, commonly referred to as biquaternions or complexified quaternions, can be expressed by the notation \cite{3}.

\parindent=0mm\vspace{0.2in}
{\it{2.1. Quaternions.}}
The quaternion algebra is an extension of the complex number to 4D algebra.
It was first invented by W. R. Hamilton in 1843 and classically denoted by
$\mathbb{H}$ in his honor. Every element of $\mathbb{H}$  has a Cartesian form given by

$$\mathbb{H}=\{q|q:=[q]_0+{\bf{i}}[q]_1+{\bf{j}}[q]_2+{\bf{k}}[q]_3,~[q]_i \in \mathbb{R},~~~i = 0, 1, 2, 3\}$$
where ${\bf{i,j, k}}$ are imaginary units obeying Hamilton’s multiplication rules (see
\cite{22})
$${\bf{i}^2}={\bf{j}^2}={\bf{k}^2}=-1$$
$${\bf{ij}}={-\bf{ji}}={\bf{k}},{\bf{jk}}={-\bf{kj}}={\bf{i}},{\bf{ki}}={-\bf{ik}}={\bf{k}}.$$

\parindent=0mm\vspace{0.2in}
Let $[q]_0$ and $q = {\bf{i}}[q]_1+{\bf{j}}[q]_2 +{\bf{k}}[q]_3 $  denote the real scalar part and the vector
part of quaternion number $q = [q]_0+{\bf{i}}[q]_1+{\bf{j}}[q]_2 +{\bf{k}}[q]_3 $, respectively. Then,
from \cite{4, upp},the real scalar part has a cyclic multiplication symmetry
$$[pql]_0=[qlp]_o=[lpq]_0,~~~\forall q,p,l \in \mathbb{H}$$
the conjugate of a quaternion $q$ is defined by $\bar{q} = [q]_0-{\bf{i}}[q]_1-{\bf{j}}[q]_2 -{\bf{k}}[q]_3 $
and the norm of $q\in \mathbb{H} $ defined as
$$|q|=\sqrt{q\bar{q}}=\sqrt{ [q]_0^2+[q]_1^2+[q]_2^2 +[q]_3^2 }.$$

\parindent=0mm\vspace{0.2in}
It is easy to verify that
$$\overline{pq}=\bar{q}\bar{p}, ~|qp|=|q||p|,~~~\forall ~q,p \in \mathbb{H}.$$

The quaternion modules $f \in L^2 (\mathbb{R}^2,\mathbb{H})$ are defined as
$$L^2 (\mathbb{R}^2,\mathbb{H}):=\{f|f:\mathbb{R}^2\rightarrow \mathbb{H},\|F\|_L^2 (\mathbb{R}^2,\mathbb{H})=\int_{\mathbb{R}^2} |f(x_1,x_2)|^2 dx_1dx_2<\infty\}.$$

\parindent=0mm\vspace{0.2in}
 An inner product of quaternion functions $f,g$ defined
on $L^2 (\mathbb{R}^2,\mathbb{H})$ given by

$$(f,g)_{L^2 (\mathbb{R}^2,\mathbb{H})}=\int_{\mathbb{R}^2}f({\bf{x}})\overline{g({\bf{x}})}d{\bf{x}}$$

\parindent=0mm\vspace{0.2in}
with symmetric real scalar part
$$\langle f,g\rangle=\frac{1}{2}\{(f,g)+(g,f)\}=\int_{\mathbb{R}^2}[f({\bf{x}})\overline{g({\bf{x}})}]_0d{\bf{x}}$$

\parindent=0mm\vspace{0.2in}
where ${\bf{x}}=(x_1,x_2)\in \mathbb{R}^2 ~ and ~ d{\bf{x}}=dx_1dx_2 .$

\parindent=0mm\vspace{0.2in}
The associated scalar norm of $f({\bf{x}})\in L^2 (\mathbb{R}^2,\mathbb{H})$ is defined by above two equations
$$\|f\|_{L^2 (\mathbb{R}^2,\mathbb{H})}^2=\langle f,f \rangle _{L^2 (\mathbb{R}^2,\mathbb{H})}=\int_{\mathbb{R}^2}|f({\bf{x}})|^2d{\bf{x}} <\infty.$$

\parindent=0mm\vspace{0.2in}
{\it{2.2. Biquaternions:}} We use the symbol $\mathbb{H_{\mathbb{C}}}$ \cite{5,6} to  denote the set of biquaternions. An  element   $\frak {h}  \in \mathbb{H_{\mathbb{C}}}$ can be expressed as
$$ \frak {h} =\frak {h}_{0} +\frak {h}_{1}\textbf{i}+\frak {h}_{2}\textbf{j}+\frak {h}_{3}\textbf{k}$$

where $\frak {h}_0,\frak {h}_1,\frak {h}_2,\frak {h}_3 \in \mathbb{C}$ are complex numbers  and $\textbf{i},\textbf{j}, \textbf{k}$ are exactly the same in real quaternions. If $\frak {h}_0=0$, then biquaternion $\frak {h}$ is known as pure biquaternion \cite{9}.

 \parindent=0mm\vspace{0.1in}
The complex numbers are generally written by the notation \cite{9} 
$$\textbf{I}^2=-1.$$
A biquaternion $q$ can also be expressed in the following way \cite{10}
$$\frak {h}=S(\frak {h})+V(\frak {h})$$
where $S(\frak {h})=\frak {h}_0$ is the complex part or scalar part of $\frak {h}$ and $V(\frak {h})=\frak {h}_{1}\textbf{i}+\frak {h}_{2}\textbf{j}+\frak {h}_{3}\textbf{k}$ is  vector part or complex quaternion-valued of $\frak {h}$.

\parindent=0mm\vspace{0.1in}
Consider two biquaternions $\frak {h}=\frak {h}_{0} +\frak {h}_{1}\textbf{i}+\frak {h}_{2}\textbf{j}+\frak {h}_{3}\textbf{k}=S(\frak {h})+V(\frak {h})$ and $\frak {g}=\frak {g}_{0} +\frak {g}_{1}\textbf{i}+\frak {g}_{2}\textbf{j}+\frak {g}_{3}\textbf{k}=S(\frak {h})+V(\frak {g}),$ then
\begin{align*}
\frak {h}+\frak {g}=[S(\frak {h})+S(\frak {g})]+[V(\frak {h})+V(\frak {g})],
\end{align*}
\begin{align*}
\lambda \frak {h}= \lambda S(\frak {h})+\lambda V(\frak {h}),
\end{align*}
\begin{align*}
\frak {h}\frak {g}=S(\frak {h})S(\frak {g})-\langle V(\frak {h}),V(\frak {g}) \rangle +S(\frak {h})V(\frak {g})+S(\frak {g})V(\frak {h})+V(\frak {h})\wedge V(\frak {g}),
\end{align*}
where $\lambda \in \mathbb{R},$
$$\langle V(\frak {h}),V(\frak {g})\rangle= \frak {h}_1\frak {g}_1+ \frak {h}_2 \frak {g}_2+\frak {h}_3\frak {g}_3,$$
$$V(\frak {h}) \wedge V(\frak {g})=(\frak {h}_2 \frak {g}_3-\frak {h}_3\frak {g}_2)\textbf{i}-(\frak {h}_1\frak {g}_3-\frak {h}_3\frak {g}_1)\textbf{j}+(\frak {h}_1\frak {g}_2-\frak {h}_2\frak {g}_1)\textbf{k}.$$
If $\frak {h}$ is orthogonal to $\frak {g}$ , then $\langle \frak {h},\frak {g} \rangle=0.$

\parindent=0mm\vspace{0.1in}
In addition, the definitions of the real and imaginary components of a biquaternion  are as follows \cite{4,8}
$$\mathcal{R}(\frak {h})=\mathcal{R}(\frak {h}_0)+\mathcal{R}(\frak {h}_1)\textbf{i}+\mathcal{R}(\frak {h}_2)\textbf{j}+\mathcal{R}(\frak {h}_3)\textbf{k},$$
$$\mathcal{I}(\frak {h})=\mathcal{I}(\frak {h}_0)+\mathcal{I}(\frak {h}_1)\textbf{i}+\mathcal{I}(\frak {h}_2)\textbf{j}+\mathcal{I}(\frak {h}_3)\textbf{k},$$

where $\mathcal{R}(\frak {h})$ and $\mathcal{I}(\frak {h})$ are real quaternions.
 $\mathcal{R}(\frak {h}_i)$ is the real part and $\mathcal{I}(\frak {h}_i)~~(i=0,1,2,3)$ is the imaginary part of a complex number. So any biquaternion   $\frak {h} \in \mathbb{H_{\mathbb{C}}}$ can be written as  $\frak {h}=\mathcal{R}(\frak {h})+\textbf{I}~\mathcal{I}(\frak {h})$. The quaternion imaginary units \textbf{i}, \textbf{j}, \textbf{k} commutes with the complex imaginary unit \textbf{I} as,  \cite{4}
 $$\textbf{i}\textbf{I}=\textbf{I}\textbf{i},\textbf{j}\textbf{I}=\textbf{I}\textbf{j},\textbf{k}\textbf{I}=\textbf{I}\textbf{k}.$$

The quaternion conjugate of a biquaternion $\frak {h} \in \mathbb{H_{\mathbb{C}}}$ is \cite{8} 
$$\widetilde{\frak {h}}=S(\frak {h})-V(\frak {h}),$$
The complex conjugate of a biquaternion $\frak {h} \in \mathbb{H_{\mathbb{C}}}$ is defined as \cite{4}
$$\check{\frak {h}}=\check{\frak {h}_0}+\check{\frak {h}_1}\textbf{i}+\check{\frak {h}_2}\textbf{j}+\check{\frak {h}_3}\textbf{k},$$

where $\check{\frak {h}_0},\check{\frak {h}_1},\check{\frak {h}_2},\check{\frak {h}_3}$ are the complex conjugates of the  coefficients of $\frak {h}$.

\parindent=0mm\vspace{0.1in}
Also, the combination of the  above two conjugations is known as biquaternion conjugate of $\frak {h}$ is defined in the following manner [26]
$$\bar{\frak {h}}=\tilde{\check{\frak {h}}}=\check{\tilde{\frak {h}}}=\check{\frak {h}_0}-\check{\frak {h}_1}\textbf{i}-\check{\frak {h}_2}\textbf{j}-\check{\frak {h}_3}\textbf{k}.$$
It is clear that the complex conjugation is multiplicative i.e,  $\check{\frak {h}\frak {g}}=\check{\frak {g}}\check{\frak {h}}$, where as  quaternion conjugation and biquaternion conjugation are involutive, that is, $\overline{\frak {g}\frak {h}}=\bar{\frak {h}}\bar{\frak {g}}$. Also, a biquaternion does not , in general, commute  with its complex conjugate \cite{19}. On reversing  the product order which give the complex conjugate result.

\parindent=0mm\vspace{0.1in}
For a  biquaternion $\frak {h}$, the norm  is defined by \cite{9}
$$||\frak {h}||^2=|\frak {h}_0|^2+|\frak {h}_1|^2+|\frak {h}_2|^2+|\frak {h}_3|^2.$$
If $||\frak {h}||=1$ , then $\frak {h}$ is said to be unit biquaternion.

\parindent=0mm\vspace{0.1in}
\hspace{0.1in} The definition of an exponential with biquaternion values for the transform's kernel is a crucial step in creating a biquaternion Fourier transform \cite{20}. In the sections that follow, we will primarily focus on exponential kernels with a biquaternion root of $-1$ \cite{8}. Numerous characteristics of the BiQFT depend on the exponential kernel \cite{4}. A biquaternion ${\gamma} \in \mathbb{H_{\mathbb{C}}}$ is a biquaternion root of  $-1$ iff $\boldsymbol{\gamma^2}=-1$ \cite{8}.

\parindent=0mm\vspace{0.1in}
For a  given  biquaternion $\frak {h}=\frak {h}_{0} +\frak {h}_{1}\textbf{i}+\frak {h}_{2}\textbf{j}+\frak {h}_{3}\textbf{k}$ and any biquaternion root $-1,\gamma$, then $\frak {h}$ can be rewritten as \cite{4}
$$\frak {h}=(\frak {h}_0'+\frak {h}_1'\gamma)+(\frak {h}_2'+\frak {h}_3'\gamma)\nu=\frak {h}_0'+\frak {h}_1'\mu+\frak {h}_2'\nu+\frak {h}_3'\gamma \nu,$$
where $\frak {h}_0', \frak {h}_1', \frak {h}_2',$ and $\frak {h}_3'$ are complex numbers and $\nu$ is a biquaternion root of $-1$ orthogonal to $\gamma$. With respect to $\gamma$ and $\nu$, the aforementioned equation enables the definition of a decomposition for every biquaternion $\frak {h}$, with $Simp(\frak {h})=(\frak {h}_0'+\frak {h}_1'\gamma)$ and $Perp(\frak {h})=(\frak {h}_2'+\frak {h}_3'\gamma)\nu$ denoting simplex and perplex parts, respectively. So 
$$\frak {h}=Simp(\frak {h})+Perp(\frak {h}).$$

The exponential of a biquaternion $\frak {h}$  is defined in the following way \cite{22}.
$$e^\frak {h}=\sum_{n\in \mathbb{N}}\frac{\frak {h}^n}{n!}.$$
A biquaternion-valued function $g(\xi) $ is expressed as \cite{4}
\begin{eqnarray*}
g(x)&=&g_0(\xi)+g_1\textbf{i}+g_2\textbf{j}+g_3\textbf{k}\\
&=& Simp(g)+Perp(g)\\
&=& S(g)+V(g)\\
&=& \mathcal{R}(g)+\textbf{I}\mathcal{I}(g),
\end{eqnarray*}

where $g_0,g_1,g_2,g_3$ are complex-valued functions. These signals can represent a wide range of physical quantities that were recorded on sensors at the same location, including magnetic rings and dipoles for electromagnetic signal recording \cite{4}.

\parindent=0mm\vspace{0.in}
\hspace{0.1in} An inner product of biquaternion functions $f,g$ can be defined on  $L^2(\mathbb{R}^2,\mathbb{H_{\mathbb{C}}})$ by
An inner product of biquaternion functions $f,g$ can be defined on  $L^2(\mathbb{R}^2,\mathbb{H_{\mathbb{C}}})$ by
$$(f,g)_{L^2(\mathbb{R}^2,\mathbb{H_{\mathbb{C}}})}=\int_{-\infty}^{\infty} \int_{-\infty}^{\infty}f(\xi_1,\xi_2)\tilde{g}(\xi_1,\xi_2)d\xi_1 d\xi_2.$$

\parindent=0mm\vspace{0.3in}
{\it{2.3 Right Sided Biquaternion Linear canonical transform}}

\parindent=0mm\vspace{0.1in}
It is well know that the biquaternion multiplication is non-commutative, as such there are three types of the biquaternion linear canonical transforms (BiQLCT); viz, the left-sided, right-sided and two-sided BiQLCT. However in the present study, we shall be focussed only on the right-sided BiQLCT \cite{20}.

\parindent=0mm\vspace{0.2in}
{\bf{Definition 2.3.}} Let $M_s=\left[ \begin{array}{cc} a_s & b_s\\ c_s & d_s \end{array} \right ] \in \mathbb{R}^{2\times 2}$ be a matrix parameter satisfying $\det|M_s|=1$, for $s=1,2$. Then RBiQLCT of a function $f \in L^2(\mathbb{R}^2,\mathbb{H}_{\mathbb{C}})$ is defined by
$$
^{\mathbb{H}_{\mathbb{C}}}\mathcal{L}_{M_1,M_2}^{RB}\{f\}(\omega,\nu)= \int_{-\infty}^{\infty}\int_{-\infty}^{\infty}f(\xi,\xi_2) K_{A_1}^{\boldsymbol{\mu}}(\xi,\omega_2)K_{A_2}^{\boldsymbol{\mu}}(\xi_2,\omega_2) d\xi_1 d\xi_2. $$ 
where the biquaternion kernels are given by

$$K_{A_i}^{\boldsymbol{\mu}}(\xi_i,\omega_i)=\left \{ \begin{array}{rcl}\frac{1}{\sqrt{2\pi|b_i|}}e^{\boldsymbol{\mu} \big (\frac{a_i}{2b_i}\xi_i^2-\frac{\xi_i\omega_i}{b_i}+\frac{d_i}{2b_i}\omega_i^2-\frac{\pi}{4}\big )}, & b_i \neq0 \\  \sqrt{d_i}e^{\boldsymbol{\mu} \frac{c_id_i}{2}\omega_i^2}\delta(\xi_i-d_i\omega_i)~~~~~~, & b_i=0
\end{array} \right.  $$

According to the definition of the BiQLCTs, we can obtain the relations between the BiQLCTs and the BiQFTs. For example, the relation between the RBiQLCT and the RBiQFT:

$$^{\mathbb{H}_{\mathbb{C}}}\mathcal{L}_{M_1,M_2}^{RB}\{f\}(\omega,\nu)=\frac{1}{2\pi \sqrt{|b_1b_2|}}F_B^R \bigg ( f(\xi_1,\xi_2)\boldsymbol{e}^{\boldsymbol{\mu}\big( \frac{a_1}{2b_1}\xi_1^2+\frac{a_2}{2b_2}\xi_2^2 \big )}\bigg )\bigg( \frac{\omega}{b_1}, \frac{\nu}{b_2} \bigg) \boldsymbol{e}^{\boldsymbol{\mu}\big( \frac{d_1}{2b_1}\omega^2+\frac{d_2}{2b_2}\nu^2-\frac{\pi}{2} \big )},$$

where the definition of the RBiQFT is the following formula:

$$F_B^R(f)(\omega,\nu)=\int_{-\infty}^{\infty}\int_{-\infty}^{\infty}f(\xi_1,\xi_2)\boldsymbol{e}^{-\boldsymbol{\mu}(\omega \xi_1+\nu \xi_2)}d\xi_1 d\xi_2.$$

\parindent=0mm\vspace{0.3in}
 {\bf{3.   Biquaternion Windowed Linear Canonical  Transform  BiQWLCT}}

\parindent=0mm \vspace{.1in}
In this section, we first introduce the notion of Biquaternion Windowed Linear Canonical Transform (BiQWLCT) and  then  various properties of BiQWLCT, such as linearity, shift, parity, orthogonality relation, inversion formula, Plancherel theorem are established.

\parindent=0mm \vspace{.1in}
{\bf{Definition 3.1.}} Let $M_s=\left[ \begin{array}{cc}
a_s& b_s\\ c_s & d_s \end{array} \right]\in \mathbb{R}^{2\times 2}$ be a matrix parameter satisfying $|M_s|=1$, for s=1,2. Let $\phi \in L^2(\mathbb{R}^2,\mathbb{H}_{\mathbb{C}})$, be a non-zero biquaternion window function   , then for any $f \in L^2(\mathbb{R}^2,\mathbb{H}_{\mathbb{C}})$, the right sided biquaternion windowed linear canonical  transform with respect to $\phi$
is defined by 
\setcounter{equation}{0}
\renewcommand{\theequation}{3.\arabic{equation}}

\begin{equation}
\mathcal{G}_{\phi}^{M_1,M_2}\{f(\boldsymbol{\xi})\}( \boldsymbol{\omega},\boldsymbol{\nu}) = \int_{-\infty}^{\infty}\int_{-\infty}^{\infty} f(\boldsymbol{\xi}) \overline{\phi(\boldsymbol{\xi} - \boldsymbol{\nu})} K_{A_1}^{\boldsymbol{\mu}}(\xi_1,\omega_1) K_{A_2}^{\boldsymbol{\mu}}(\xi_2,\omega_2) \, d\xi_1 \, d\xi_2
\label{eq:3.1}
\end{equation}

where $(\boldsymbol{\omega},\boldsymbol{\nu})\in \mathbb{R}^2$ and the kernel functions are defined by
$$K_{M_1}^{\boldsymbol{\mu}}(\xi_1,\omega_1)=\left \{ \begin{array}{rcl}\frac{1}{\sqrt{2\pi|b_1|}}e^{\boldsymbol{\mu} \big (\frac{a_1}{2b_1}\xi_1^2-\frac{\xi_1\omega_1}{b_1}+\frac{d_1}{2b_1}\omega_1^2-\frac{\pi}{4}\big )}, & b_1 \neq0 \\  \sqrt{d_1}e^{\boldsymbol{\mu}\frac{c_1d_1}{2}\omega_1^2}\delta(\xi_1-d_1\omega_1)~~~~~~, & b_1=0
\end{array} \right.$$
and
$$K_{M_2}^{\boldsymbol{\vartheta}}(\xi_2,\omega_2)=\left \{ \begin{array}{rcl}\frac{1}{\sqrt{2\pi|b_2|}}e^{\boldsymbol{\vartheta} \big (\frac{a_2}{2b_2}\xi_2^2-\frac{\xi_2\omega_2}{b_2}+\frac{d_2}{2b_2}\omega_2^2-\frac{\pi}{4}\big )}, & b_2 \neq0 \\  \sqrt{d_2}e^{\boldsymbol{\vartheta}\frac{c_2d_2}{2}\omega_2^2}\delta(\xi_2-d_2\omega_2)~~~~~~, & b_2=0
\end{array} \right.$$
where $\delta(\boldsymbol{\xi})$ representing the Dirac function.

\parindent=0mm \vspace{.1in} The Dirac  delta can be losely thought of as a function on the real line which is zero everywhere except at the origin, where it is finite,
$$\delta(\xi)= \left \{ \begin{array}{rcl} +\infty ~~~~, & \xi=0 \\  0~~~~~ , & \xi \neq 0
\end{array} \right.$$
and which is also constrained to satissfy the identity 
$$\int_\infty^\infty \delta(\xi)d\xi=1. $$

In the above defintion, for $b_i=0$, $i=1,2$, the right sided BiQWLCT of a signal is essentially a chirp multiplication and it is of no  interest  in this paper. Therefore in the rest of the paper, we assume $b_i \ne 0$. Under some suitable conditions, the right sided BiQWLCT above is invertible and the inverse is given as:

\parindent=0mm \vspace{.1in} For a fixed $\boldsymbol{\nu}$, we have
$$\mathcal{G}_{\phi}^{M_1,M_2}\{f(\boldsymbol{\xi})\}( \boldsymbol{\omega},\boldsymbol{\nu})=~^{\mathbb{H}_{\mathbb{C}}}\mathcal{L}_{M_1,M_2}^{RB}\{f(\boldsymbol{\xi})\overline{\phi(\boldsymbol{\xi}-\boldsymbol{\nu})}\}(\boldsymbol{\omega}).$$
The inverse BiQWLCT is given as:
\begin{align}
f(\boldsymbol{\xi})\overline{\phi(\boldsymbol{\xi}-\boldsymbol{\nu})}&= ~^{\mathbb{H}_{\mathbb{C}}}\mathcal{L}_{A_1^{-1},A_2^{-1}}^{RB} \{ \mathcal{G}_{\phi}^{M_1,M_2}[f(\boldsymbol{\xi})]( \boldsymbol{\omega},\boldsymbol{\nu})\}\\
&=\int_{-\infty}^{\infty}\int_{-\infty}^{\infty}\mathcal{G}_{\phi}^{M_1,M_2}\{f(\boldsymbol{\xi})\}( \boldsymbol{\omega},\boldsymbol{\nu})\overline{K_{A_1^{-1}}^{\boldsymbol{\mu}}(\xi_1,\omega_1)K_{A_2^{-1}}^{\boldsymbol{\mu}}(\xi_2,\omega_2)} d\omega_1\omega_2
\end{align}

\parindent=0mm \vspace{.1in}
{\bf{ 3.2. Properties of right sided BiQWLCT}}

\parindent=0mm \vspace{.2in}
In this subsection, we discuss several basic properties of the right sided BiQWLCT. These properties play important roles in signal representation.

\parindent=0mm \vspace{.1in}
{\bf{Theorem 3.2.}} ( Linearity) Let $\phi \in L^2(\mathbb{R}^2,\mathbb{H}_{\mathbb{C}})$, be a biquaternion window function  and $f,g \in L^2(\mathbb{R}^2,\mathbb{H}_{\mathbb{C}})$.Then the right sided BiQWLCT is a linear operator, namely, 
$$\mathcal{G}_{\phi}^{M_1,M_2}\{\alpha f(\boldsymbol{\xi})+\beta g(\boldsymbol{\xi})\}( \boldsymbol{\omega},\boldsymbol{\nu})=\alpha \mathcal{G}_{\phi}^{M_1,M_2}\{f\}( \boldsymbol{\omega},\boldsymbol{\nu})+\beta \mathcal{G}_{\phi}^{M_1,M_2} \{g\}( \boldsymbol{\omega},\boldsymbol{\nu})$$ 
where $\alpha,\beta$ are arbitrary real constants. This property sheds lights on the right sided BiQWLCT analysis of multi-component signals.

\parindent=0mm \vspace{.2in}
\textit{Proof}.  By using the linearity of integration involved in $\ref{eq:3.1}$, the proof follows .

\parindent=0mm \vspace{.1in}
{\bf{Theorem 3.3.}} (Inversion formula) Let $f \in L^2(\mathbb{R}^2,\mathbb{H}_{\mathbb{C}})$ and $\phi \in L^2(\mathbb{R}^2,\mathbb{H}_{\mathbb{C}})$, be a biquaternion window function ,$0<||\phi||^2<\infty$. Then we have the inversion formula of BiQWLCT ,
\begin{eqnarray*}
f\boldsymbol{(\xi)}&=&\bigg(\mathcal{G}_{\phi}^{M_1,M_2}\{f\}( \boldsymbol{\omega},\boldsymbol{\nu})\bigg)^{-1}\bigg[\mathcal{G}_{\phi}^{M_1,M_2}\{f\}( \boldsymbol{\omega},\boldsymbol{\nu})\bigg] \\
&=& \frac{1}{||\phi||^2}\int_{-\infty}^{\infty}\int_{-\infty}^{\infty}\int_{-\infty}^{\infty}\int_{-\infty}^{\infty}\mathcal{G}_{\phi}^{M_1,M_2}\{f(\boldsymbol{\xi})\}( \boldsymbol{\omega},\boldsymbol{\nu})\overline{K_{M_1^{-1}}^{\boldsymbol{\mu}}(\xi_1,\omega_1)}~\overline{K_{M_2^{-1}}^{\boldsymbol{\mu}}(\xi_2,\omega_2)} d\omega_1\omega_2 d\nu_1 d\nu_2.
\end{eqnarray*}
where $M_s^{-1}=\left[ \begin{array}{cc}
d_s& -b_s\\ -c_s & a_s \end{array} \right]\in \mathbb{R}^{2\times 2}$, for $s=1,2$.

\parindent=0mm \vspace{.1in}
\textit{Proof}. Based on the definition of right sided BiQWLCT, multiply $(3.3)$ from the right side by $\phi(\boldsymbol{\xi}-\boldsymbol{\nu})$ and integrating with respect to $d\boldsymbol{\nu}$ , we get
\begin{align*}
\int_{-\infty}^{\infty}\int_{-\infty}^{\infty}f(\boldsymbol{\xi})\overline{\phi(\boldsymbol{\xi}-\boldsymbol{\nu})}\phi(\boldsymbol{\xi}-\boldsymbol{\nu})d{\nu}_1d{\nu}_2 &= \int_{-\infty}^{\infty}\int_{-\infty}^{\infty}f(\boldsymbol{\xi})|\phi(\boldsymbol{\xi}-\boldsymbol{\nu})|^2d{\nu}_1d{\nu}_2\\
&= \int_{-\infty}^{\infty}\int_{-\infty}^{\infty}\int_{-\infty}^{\infty}\int_{-\infty}^{\infty}\mathcal{G}_{\phi}^{M_1,M_2}\{f(\boldsymbol{\xi})\}( \boldsymbol{\omega},\boldsymbol{\nu})\phi(\boldsymbol{\xi}-\boldsymbol{\nu})\\
&\qquad\qquad \times \overline{K_{M_1^{-1}}^{\boldsymbol{\mu}}(\xi_1,\omega_1)}~\overline{K_{M_2^{-1}}^{\boldsymbol{\mu}}(\xi_2,\omega_2)} d\omega_1\omega_2 d\nu_1 d\nu_2.
\end{align*}

The associated scalar norm of $f(\boldsymbol{\xi}) \in L^2(\mathbb{R}^2,\mathbb{H}_{\mathbb{C}})$ is given by 
$$||f||_{L^2(\mathbb{R}^2,\mathbb{H}_{\mathbb{C}})}^2=\langle f,f \rangle_{L^2(\mathbb{R}^2,\mathbb{H}_{\mathbb{C}})}=\int_{-\infty}^{\infty}\int_{-\infty}^{\infty}|f(\boldsymbol{\xi})|^2d\boldsymbol{\xi} <\infty.$$
We can write $\int_{-\infty}^{\infty}\int_{-\infty}^{\infty}|\phi(\boldsymbol{\xi}-\boldsymbol{\nu})|^2d{\nu}_1d{\nu}_2=||\phi||^2$
 and using the above result , we have
\begin{eqnarray*}
f\boldsymbol{(\xi)}
= \frac{1}{||\phi||^2}\int_{-\infty}^{\infty}\int_{-\infty}^{\infty}\int_{-\infty}^{\infty}\int_{-\infty}^{\infty}\mathcal{G}_{\phi}^{M_1,M_2}\{f(\boldsymbol{\xi})\}( \boldsymbol{\omega},\boldsymbol{\nu})\overline{K_{M_1^{-1}}^{\boldsymbol{\mu}}(\xi_1,\omega_1)}~\overline{K_{M_2^{-1}}^{\boldsymbol{\mu}}(\xi_2,\omega_2)} d\omega_1\omega_2 d\nu_1 d\nu_2.
\end{eqnarray*}

which completes the proof.

\parindent=0mm \vspace{.1in}
{\bf{Theorem 3.4.}} (Plancherel theorem) Let $\phi \in L^2(\mathbb{R}^2,\mathbb{H}_{\mathbb{C}})$, be a biquaternion window function  and let $f,g \in L^2(\mathbb{R}^2,\mathbb{H}_{\mathbb{C}})$. Then we have
$$\big \langle \mathcal{G}_{\phi}^{M_1,M_2}\{f\}( \boldsymbol{\omega},\boldsymbol{\nu}), \mathcal{G}_{\phi}^{M_1,M_2}\{g\}( \boldsymbol{\omega},\boldsymbol{\nu})\rangle _{L^2(\mathbb{R}^2,\mathbb{H}_{\mathbb{C}})}=\langle f,g \big \rangle_{L^2(\mathbb{R}^2,\mathbb{H}_{\mathbb{C}})}.$$

\parindent=0mm \vspace{.1in}
\textit{Proof}. Applying the inversion formula and the definition of the right sided BiQWLCT, we obtain 
\begin{align*}
\left\langle \mathcal{G}_{\phi}^{M_1,M_2}\{f\}( \boldsymbol{\omega},\boldsymbol{\nu}), \mathcal{G}_{\phi}^{M_1,M_2}\{g\}( \boldsymbol{\omega},\boldsymbol{\nu})\right\rangle &= \displaystyle\int_{-\infty}^{\infty}\int_{-\infty}^{\infty}\mathcal{G}_{\phi}^{M_1,M_2}\{f\}( \boldsymbol{\omega},\boldsymbol{\nu})\overline{\mathcal{G}_{\phi}^{M_1,M_2}\{g\}}( \boldsymbol{\omega},\boldsymbol{\nu})d\omega_1 d\omega_2\\\\
& = \displaystyle \int_{-\infty}^{\infty}\int_{-\infty}^{\infty}\int_{-\infty}^{\infty}\int_{-\infty}^{\infty}\mathcal{G}_{\phi}^{M_1,M_2}\{f(\boldsymbol{\xi})\}( \boldsymbol{\omega},\boldsymbol{\nu})\mathcal{G}_{\phi}^{M_1,M_2}\{g(\boldsymbol{\xi})\}( \boldsymbol{\omega},\boldsymbol{\nu})\\\\
&\qquad\quad\qquad \times  \phi(\boldsymbol{\xi}-\boldsymbol{\nu}) \overline{K_{M_1^{-1}}^{\boldsymbol{\mu}}(\xi_1,\omega_1)}~\overline{K_{M_2^{-1}}^{\boldsymbol{\mu}}(\xi_2,\omega_2)} d\omega_1\omega_2 d\xi_1 d\xi_2\\\\
& =  \displaystyle \int_{-\infty}^{\infty}\int_{-\infty}^{\infty}\bigg \langle \mathcal{G}_{\phi}^{M_1,M_2}\{f(\boldsymbol{\xi})\}( \boldsymbol{\omega},\boldsymbol{\nu})\int_{-\infty}^{\infty}\int_{-\infty}^{\infty}\mathcal{G}_{\phi}^{M_1,M_2}\{g(\boldsymbol{\xi})\}( \boldsymbol{\omega},\boldsymbol{\nu})\\\\
&\qquad\qquad\quad \times  \phi(\boldsymbol{\xi}-\boldsymbol{\nu}) \overline{K_{M_1^{-1}}^{\boldsymbol{\mu}}(\xi_1,\omega_1)}~\overline{K_{M_2^{-1}}^{\boldsymbol{\mu}}(\xi_2,\omega_2)}  d\xi_1 d\xi_2 \bigg \rangle d\omega_1\omega_2\\\\
& = \displaystyle \int_{-\infty}^{\infty}\int_{-\infty}^{\infty}f(\boldsymbol{\xi})\bar{g}(\boldsymbol{\xi})d\boldsymbol{\xi}\\\\
&= \langle f,g \rangle.
\end{align*}

which completes the proof.

\parindent=0mm \vspace{.1in}
{\bf{Corollary 3.5:}} If $f=g$ then
$$\int_{-\infty}^{\infty}\int_{-\infty}^{\infty}\bigg | \mathcal{G}_{\phi}^{M_1,M_2}\{f\}( \boldsymbol{\omega},\boldsymbol{\nu})\bigg |^2d\omega_1d\nu_1=||f||_{ L^2(\mathbb{R}^2,\mathbb{H}_{\mathbb{C}})}^2.$$

\parindent=0mm \vspace{.1in}
{\bf{Theorem 3.6.}}(Parity) Let $\phi \in L^2(\mathbb{R}^2,\mathbb{H}_{\mathbb{C}})$, be a biquaternion window function  and let $f \in L^2(\mathbb{R}^2,\mathbb{H}_{\mathbb{C}})$. Then we have

$$\mathcal{G}_{\mathcal{P}{\phi}}^{M_1,M_2}\{\mathcal{P}f\}( \boldsymbol{\omega},\boldsymbol{\nu})=\mathcal{G}_{\mathcal{P}{\psi}}^{M_1,M_2}( -\boldsymbol{\omega},-\boldsymbol{\nu})$$
where $Pf=f(-\boldsymbol{\xi})$ and $\mathcal{P}\phi=\phi(-\boldsymbol{\xi})$ for every window function $\phi \in L^2(\mathbb{R}^2,\mathbb{H}_{\mathbb{C}}) $.

\parindent=0mm \vspace{.1in}
\textit{Proof}. For every $f \in L^2(\mathbb{R}^2,\mathbb{H}_{\mathbb{C}})$, we have
\begin{align*}
\mathcal{G}_{\mathcal{P}{\phi}}^{M_1,M_2}\{\mathcal{P}f\}( \boldsymbol{\omega},\boldsymbol{\nu})&=\int_{-\infty}^{\infty}\int_{-\infty}^{\infty}f(\boldsymbol{-\xi})\overline{\phi(-(\boldsymbol{\xi}-\boldsymbol{\nu}))}K_{M_1}^{\boldsymbol{\mu}}(\xi_1,\omega_1) K_{M_2}^{\boldsymbol{\mu}}(\xi_2,\nu_1)d\xi_1d\xi_2\\\\
&=\int_{-\infty}^{\infty}\int_{-\infty}^{\infty}f(\boldsymbol{-\xi})\overline{\phi(-(\boldsymbol{\xi}-\boldsymbol{(-\nu)}))}\frac{1}{\sqrt{2\pi b_1}}e^{\boldsymbol{\mu}\big (\frac{a_1}{2b_1}(-\xi_1)^2-\frac{(-\xi_1)(-\omega_1)}{b_1}+\frac{d_1}{2b_1}(-\omega_1)^2-\frac{\pi}{4}\big )} \\ \\
&\qquad\qquad\qquad\qquad\qquad\quad \times \frac{1}{\sqrt{2\pi b_2}} e^{\boldsymbol{\mu}\big (\frac{a_2}{2b_2}(-\xi_2)^2-\frac{(-\xi_2)(-\nu_1)}{b_2}+\frac{d_2}{2b_2}(-\nu_1)^2-\frac{\pi}{4}\big )} d\xi_1d\xi_2\\ \\
&= \mathcal{G}_{\mathcal{P}{\psi}}^{M_1,M_2}( -\boldsymbol{\omega},-\boldsymbol{\nu}).
\end{align*}
Hence, we get the desired result.

\parindent=0mm \vspace{.1in}
{\bf{Theorem 3.7.}} (Orthogonality relation) Let $f,g  \in L^2(\mathbb{R}^2,\mathbb{H}_{\mathbb{C}})$ and $\phi, \Psi   \in L^2(\mathbb{R}^2,\mathbb{H}_{\mathbb{C}})$ be a non-zero biquaternion window function, then we have
$$\int_{-\infty}^{\infty}\int_{-\infty}^{\infty}\int_{-\infty}^{\infty}\int_{-\infty}^{\infty} \bigg[ \mathcal{G}_{\phi}^{M_1,M_2}\{f\}( \boldsymbol{\omega},\boldsymbol{\nu})\overline{\mathcal{G}_{\Psi}^{M_1,M_2}\{g\}( \boldsymbol{\omega},\boldsymbol{\nu})}\bigg]_0 d\omega_1 d\omega_2 d\nu_1 d\nu_2= \big [\langle \Upsilon_{\phi, \Psi}f,g\rangle \big ]_0.$$

\textit{Proof}. From the definition of RBiQWLCT
\begin{align*}
& \int_{-\infty}^{\infty}\int_{-\infty}^{\infty}\int_{-\infty}^{\infty}\int_{-\infty}^{\infty} \bigg[ \mathcal{G}_{\phi}^{M_1,M_2}\{f\}( \boldsymbol{\omega},\boldsymbol{\nu})\overline{\mathcal{G}_{\Psi}^{M_1,M_2}\{g\}( \boldsymbol{\omega},\boldsymbol{\nu})}\bigg]_0 d\omega_1 d\omega_2 d\nu_1 d\nu_2\\\\
&= \int_{-\infty}^{\infty}\int_{-\infty}^{\infty}\int_{-\infty}^{\infty}\int_{-\infty}^{\infty} \bigg[ \mathcal{G}_{\phi}^{M_1,M_2}\{f\}( \boldsymbol{\omega},\boldsymbol{\nu})\\\\
&\qquad\qquad\qquad\times \int_{-\infty}^{\infty}\int_{-\infty}^{\infty} \overline{g(\boldsymbol{\xi})}~ \overline{\Psi(\boldsymbol{\xi}-\boldsymbol{\nu})}\frac{1}{\sqrt{2\pi|b_1|}}e^{-\boldsymbol{\mu} \big (\frac{a_1}{2b_1}\xi_1^2-\frac{\xi_1\omega_1}{b_1}+\frac{d_1}{2b_1}\omega_1^2-\frac{\pi}{4}\big )}\\\\
&\qquad\qquad\qquad\qquad\qquad\qquad\times \frac{1}{\sqrt{2\pi|b_1|}}e^{-\boldsymbol{\mu} \big (\frac{a_2}{2b_2}\xi_2^2-\frac{\xi_2\omega_2}{b_2}+\frac{d_2}{2b_2}\omega_2^2-\frac{\pi}{4}\big )}\bigg ]_0 d\omega_1 d\omega_2 d\nu_1 d\nu_2 d\xi_1 d\xi_2\\\\
&= \int_{-\infty}^{\infty}\int_{-\infty}^{\infty}\int_{-\infty}^{\infty}\int_{-\infty}^{\infty} \bigg[ \mathcal{G}_{\phi}^{M_1,M_2}\{f\}( \boldsymbol{\omega},\boldsymbol{\nu})\\\\
&\qquad\qquad\qquad\qquad \times \overline{\int_{-\infty}^{\infty}\int_{-\infty}^{\infty}g(\boldsymbol{\xi})\overline{\Psi(\boldsymbol{\xi}-\boldsymbol{\nu})}K_{M_1}^{\boldsymbol{\mu}}(\xi_1,\omega_1) K_{M_2}^{\boldsymbol{\mu}}(\xi_2,\omega_2) }\bigg]_0 d\omega_1 d\omega_2 d\nu_1 d\nu_2 d\xi_1 d\xi_2\\\\
\end{align*}
\begin{align*}
&= \int_{-\infty}^{\infty}\int_{-\infty}^{\infty}\int_{-\infty}^{\infty}\int_{-\infty}^{\infty} \bigg[ \mathcal{G}_{\phi}^{M_1,M_2}\{f\}( \boldsymbol{\omega},\boldsymbol{\nu})\\\\
& \qquad\qquad\qquad\qquad \times \int_{-\infty}^{\infty}\int_{-\infty}^{\infty} \overline{g(\boldsymbol{\xi})}\Psi(\boldsymbol{\xi}-\boldsymbol{\nu})K_{M_1}^{-\boldsymbol{\mu}}(\xi_1,\omega_1) K_{M_2}^{\boldsymbol{-\mu}}(\xi_2,\omega_2) \bigg]_0 d\omega_1 d\omega_2 d\nu_1 d\nu_2 d\xi_1 d\xi_2\\\\
&= \int_{-\infty}^{\infty}\int_{-\infty}^{\infty}\int_{-\infty}^{\infty}\int_{-\infty}^{\infty} \bigg[ \int_{-\infty}^{\infty}\int_{-\infty}^{\infty} \mathcal{G}_{\phi}^{M_1,M_2}\{f\}( \boldsymbol{\omega},\boldsymbol{\nu})K_{M_1}^{-\boldsymbol{\mu}}(\xi_1,\omega_1) K_{M_2}^{\boldsymbol{-\mu}}(\xi_2,\omega_2)\\\\
& \qquad\qquad\qquad\qquad \qquad\qquad\qquad\qquad \qquad\qquad\qquad\qquad \cdot \Psi(\boldsymbol{\xi}-\boldsymbol{\nu})\overline{g(\boldsymbol{\xi})}\bigg]_0 d\omega_1 d\omega_2 d\nu_1 d\nu_2 d\xi_1 d\xi_2.
\end{align*}
It is know that $K_{M_1}^{\boldsymbol{\mu}}(\xi_1,\omega_1), K_{M_2}^{\boldsymbol{\mu}}(\xi_2,\omega_2)  $ satisfies the conditions:
\begin{itemize}
    \item \( \overline{K_{M_1}^{\boldsymbol{\mu}}(\xi_1,\omega_1) } = K_{M_1}^{-\boldsymbol{\mu}}(\xi_1,\omega_1) \)
    \item \( K_{M_1}^{-\boldsymbol{\mu}}(\xi_1,\omega_1) = K_{M_1^{-1}}^{\boldsymbol{\mu}}(\xi_1,\omega_1) \)
\end{itemize}

\begin{itemize}
    \item \( \overline{K_{M_2}^{\boldsymbol{\mu}}(\xi_2,\omega_2) } = K_{M_2}^{-\boldsymbol{\mu}}(\xi_2,\omega_2) \)
    \item \( K_{M_2}^{-\boldsymbol{\mu}}(\xi_2,\omega_2) = K_{M_2^{-1}}^{\boldsymbol{\mu}}(\xi_2,\omega_2) \)
\end{itemize}
 Using inversion formula and the above kernel properties, we obtain
 \begin{align*}
 & \int_{-\infty}^{\infty}\int_{-\infty}^{\infty}\int_{-\infty}^{\infty}\int_{-\infty}^{\infty} \bigg[ \mathcal{G}_{\phi}^{M_1,M_2}\{f\}( \boldsymbol{\omega},\boldsymbol{\nu})\overline{\mathcal{G}_{\Psi}^{M_1,M_2}\{g\}( \boldsymbol{\omega},\boldsymbol{\nu})}\bigg]_0 d\omega_1 d\omega_2 d\nu_1 d\nu_2\\\\
 &=\int_{-\infty}^{\infty}\int_{-\infty}^{\infty}\int_{-\infty}^{\infty}\int_{-\infty}^{\infty} \bigg[ \int_{-\infty}^{\infty}\int_{-\infty}^{\infty} \mathcal{G}_{\phi}^{M_1,M_2}\{f\}( \boldsymbol{\omega},\boldsymbol{\nu})K_{M_1^{-1}}^{\boldsymbol{\mu}}(\xi_1,\omega_1) K_{M_2^{-1}}^{\boldsymbol{\mu}}(\xi_2,\omega_2)d\omega_1 d\omega_2 \\\\
& \qquad\qquad\qquad\qquad \qquad\qquad\qquad\qquad \qquad\qquad\qquad\qquad \times \Psi(\boldsymbol{\xi}-\boldsymbol{\nu})\overline{g(\boldsymbol{\xi})}\bigg]_0 d\nu_1 d\nu_2 d\xi_1 d\xi_2\\\\
&=\int_{-\infty}^{\infty}\int_{-\infty}^{\infty}\int_{-\infty}^{\infty}\int_{-\infty}^{\infty} \big[ f(\boldsymbol{\xi}) \overline{g(\boldsymbol{\xi})\Psi(\boldsymbol{\xi}-\boldsymbol{\nu})}\Psi(\boldsymbol{\xi}-\boldsymbol{\nu})\big]_0 d\xi_1 d\xi_2 d\nu_1 d\nu_2\\\\
&= \bigg[ \int_{-\infty}^{\infty}\int_{-\infty}^{\infty}f(\boldsymbol{\xi})\overline{g(\boldsymbol{\xi})}d\xi_1 d\xi_2 \int_{-\infty}^{\infty}\int_{-\infty}^{\infty}|\Psi(\boldsymbol{\xi}-\boldsymbol{\nu})|^2 d\nu_1 d\nu_2\bigg]_0\\\\
&= \bigg [\int_{-\infty}^{\infty}\int_{-\infty}^{\infty}\Upsilon_{\phi,\Psi}f(\boldsymbol{\xi})\overline{g(\boldsymbol{\xi})}d\xi_1 d\xi_2\bigg]_0\\\\
&= \big [ \big \langle \Upsilon_{\phi,\Psi}f,g \big \rangle \big ]_0,
 \end{align*}
where $\Upsilon_{\phi,\Psi}=\displaystyle\int_{-\infty}^{\infty}\int_{-\infty}^{\infty}|\Psi(\boldsymbol{\xi}-\boldsymbol{\nu})|^2 d\nu_1 d\nu_2$.
Thus the proof is completed.

\parindent=0mm \vspace{.1in}
{\bf{Theorem 3.8.}}(Shift)  Let $\phi \in L^2(\mathbb{R}^2,\mathbb{H}_{\mathbb{C}})$, be a biquaternion window function  and let $f \in L^2(\mathbb{R}^2,\mathbb{H}_{\mathbb{C}})$. Assume that $m_s=\omega_s-a_sr_s,~~n_s=\nu_s-r_s,~~s=1,2$. then we obtain
$$\mathcal{G}_{\phi}^{M_1,M_2}\{f(\boldsymbol{\xi}-\boldsymbol{r})\}( \boldsymbol{\omega},\boldsymbol{\nu})=\mathcal{G}_{\phi}^{M_1,M_2}\{f\}( \boldsymbol{m},\boldsymbol{n})e^{\boldsymbol{\mu}r_1\omega_1c_1}e^{\boldsymbol{\mu}r_2\omega_2c_2}e^{\boldsymbol{-\mu}\frac{a_1r_1^2}{2}c_1}e^{\boldsymbol{-\mu}\frac{a_2r_2^2}{2}c_2},$$
where $  \textbf{r}=(r_1,r_2),~ \textbf{m}=(m_1,m_2)$, and $ \textbf{n}=(n_1,n_2) \in \mathbb{R}^2$.

\parindent=0mm \vspace{.1in}
\textit{Proof}. From the definition $\ref{eq:3.1}$, we have
\begin{eqnarray*}
\mathcal{G}_{\phi}^{M_1,M_2}\{f(\boldsymbol{\xi}-\boldsymbol{r})\}( \boldsymbol{\omega},\boldsymbol{\nu})=\int_{-\infty}^{\infty}\int_{-\infty}^{\infty}f(\boldsymbol{\xi})\overline{\phi(\boldsymbol{\xi}-\boldsymbol{\nu})}K_{A_1}^{\boldsymbol{\mu}}(\xi_1,\omega_1) K_{A_2}^{\boldsymbol{\mu}}(\xi_2,\omega_2)d\xi_1d\xi_2.
\end{eqnarray*}
By making the change of variable $\boldsymbol{\xi}=\textbf{t}-\textbf{r}$
 in the above expression, we obtain
\begin{align*}
\mathcal{G}_{\phi}^{M_1,M_2}\{f(\boldsymbol{\xi}-\boldsymbol{r})\}( \boldsymbol{\omega},\boldsymbol{\nu})&=\int_{-\infty}^{\infty}\int_{-\infty}^{\infty} f(\textbf{t}-\textbf{r})\overline{\psi(\textbf{t}-(\boldsymbol{\nu}-\boldsymbol{r}))}\frac{1}{\sqrt{2\pi b_1}}e^{\boldsymbol{\mu}\big ( \frac{a_1}{2b_1}(t_1-r_1)^2-\frac{(t_1-r_1)}{b_1}\omega_1+\frac{d_1}{2b_1}\omega_1^2-\frac{\pi}{4}\big )}\\\\
&\qquad\qquad\qquad\qquad\qquad\qquad \times  \frac{1}{\sqrt{2\pi b_1}}e^{\boldsymbol{\mu}\big ( \frac{a_2}{2b_2}(t_2-r_2)^2-\frac{(t_2-r_2)}{b_2}\nu_1+\frac{d_2}{2b_2}\nu_1^2-\frac{\pi}{4}\big )}dt_1dt_2\\\\
&= \int_{-\infty}^{\infty}\int_{-\infty}^{\infty} f(\textbf{t}-\textbf{r})\overline{\psi(\textbf{t}-(\boldsymbol{\nu}-\boldsymbol{r}))}\frac{1}{\sqrt{2\pi b_1}}e^{\boldsymbol{\mu}\big( \frac{a_1}{2b_1}r_1^2-\frac{r_1\omega_1}{b_1}\big )}\cdot e^{\boldsymbol{\mu}\big( \frac{a_1}{2b_1}r_1^2-\frac{r_1\omega_1}{b_1}\big )}\\\\
&\qquad\qquad\qquad\qquad\qquad \times  \frac{1}{\sqrt{2\pi b_1}} e^{\boldsymbol{\mu}\big(\frac{a_1}{2b_1}t_1^2-\frac{(\omega_1-a_1r_1)}{b_1}t_1+\frac{d_1}{2b_1}(\omega_1-r_1a_1)^2-\frac{\pi}{4}\big )} dt_1 dt_2\\\\
&\qquad\qquad\qquad\qquad\qquad\qquad\qquad\qquad \times   e^{\boldsymbol{\mu}\big(\frac{a_2}{2b_2}t_2^2-\frac{(\omega_2-a_2r_2)}{b_2}t_2+\frac{d_2}{2b_2}(\omega_2-r_2a_2)^2-\frac{\pi}{4}\big )}
\end{align*}

Applying the definition of the right sided BiQWLCT, we finally arrive at
$$\mathcal{G}_{\phi}^{M_1,M_2}\{f(\boldsymbol{\xi}-\boldsymbol{r})\}( \boldsymbol{\omega},\boldsymbol{\nu})=\mathcal{G}_{\phi}^{M_1,M_2}\{f\}( \boldsymbol{m},\boldsymbol{n})e^{\boldsymbol{\mu}r_1\omega_1c_1}e^{\boldsymbol{\mu}r_2\omega_2c_2}e^{\boldsymbol{-\mu}\frac{a_1r_1^2}{2}c_1}e^{\boldsymbol{-\mu}\frac{a_2r_2^2}{2}c_2}$$
which  completes the proof.

\parindent=0mm \vspace{.2in}
{\bf{4. Heisenberg Uncertainty Principle for the Biquaternion Windowed Linear Canonical  Transform}}

\parindent=0mm \vspace{.1in}
Quaternions and the Heisenberg uncertainty principle are fundamental concepts in quantum physics. According to the uncertainty principle of quantum physics, it is impossible to know an electron's position and velocity at the same time. In other words, when one learns more about a particle's position, one learns less about its momentum or velocity. According to the uncertainty principle in signal processing, the sum of the signal's variances in the time and frequency domains has a lower bound.In Refs. \cite{27}, the Fourier transform's and the quaternion Fourier transform's uncertainty concepts were explored. In Refs. \cite{28,29,30, ap}, the authors recently established the uncertainty principles related to the LCT.
Reference \cite{25} provided a discussion of the WLCT's uncertainty principles. They used Lieb's uncertainty principles as a generalisation for the WLCT domains. Recent research has expanded Heisenberg's uncertainty relations to the Biquaternion linear canonical transform \cite{31}. This uncertainty principle sets a lower bound on the sum of the effective widths of quaternion-valued signals in the spatial and frequency domains. The Heisenberg uncertainty principles of the QWLCT have been studied in \cite{kk}. Based on the Heisenberg uncertainty principles of the  QWLCT, we can obtain the Heisenberg uncertainty principle of the right sided BiQWLCT

\parindent=0mm \vspace{.1in}
{\bf{Lemma 4.1.}} Let $M_s=\left[ \begin{array}{cc}
a_s& b_s\\ c_s & d_s \end{array} \right]\in \mathbb{R}^{2\times 2}$ be a matrix parameter satisfying $|M_s|=1$, for s=1,2. Let $f \in L^2(\mathbb{R}^2,\mathbb{H}_{\mathbb{C}})$ and 
$\omega_k ~^{\mathbb{H}_{\mathbb{C}}}\mathcal{L}_{M_1,M_2}^{RB}\{f\} \in L^2(\mathbb{R}^2,\mathbb{H}_{\mathbb{C}}) $, then we obtain

\setcounter{equation}{0}
\renewcommand{\theequation}{4.\arabic{equation}}

\setcounter{equation}{0}
\renewcommand{\theequation}{4.\arabic{equation}}

\begin{align}
& \int_{-\infty}^{\infty}\int_{-\infty}^{\infty}\xi_k^2|f(\xi_1,\xi_2)|^2 d \xi_1 d \xi_2 \nonumber 
 \int_{-\infty}^{\infty}\int_{-\infty}^{\infty} \omega_k^2|\, ^{\mathbb{H}_{\mathbb{C}}}\mathcal{L}_{M_1,M_2}^{RB}\{f\}(\omega,\nu)|^2 d\omega d\nu \nonumber \\
& \quad\quad\quad\quad\quad\quad\quad\quad\quad\quad\quad\quad\quad\quad\quad\quad\quad\quad\quad \geq \bigg \{ \frac{b_k}{2}\int_{-\infty}^{\infty}\int_{-\infty}^{\infty}|f(\xi_1,\xi_2)|^2 d \xi_1 d \xi_2 \bigg \}^2 \label{eq:4.1}
\end{align}

where $k=1,2$.

\parindent=0mm \vspace{.1in}
The above equality holds if and only if $f$ is the Gaussian function, namely,
$$f(\xi_1,\xi_2)=C_0e^{-(\alpha_1\xi_1^2+\alpha_2\xi_2^2)}$$
where $C_0$ is a complex constant and $\alpha_1, \alpha_2 \in \mathbb{R}$ are positive real constants.

\parindent=8mm \vspace{.1in}
Substituting the inverse transform for the right sided BiQLCT into the left hand side of \ref{eq:4.1}, we obtain

\begin{align*}
& \int_{-\infty}^{\infty}\int_{-\infty}^{\infty}\xi_k^2|~^{\mathbb{H}_{\mathbb{C}}}\mathcal{L}_{M_1,M_2}^{-1} [\, ^{\mathbb{H}_{\mathbb{C}}}\mathcal{L}_{M_1,M_2}^{RB}\{f\}](\xi_1,\xi_2)|^2 d \xi_1 d \xi_2 \nonumber 
 \int_{-\infty}^{\infty}\int_{-\infty}^{\infty} \omega_k^2|\,  ^{\mathbb{H}_{\mathbb{C}}}\mathcal{L}_{M_1,M_2}^{RB}\{f\}(\omega,\nu)|^2 d\omega d\nu \nonumber \\
& \quad\quad\quad\quad\quad\quad\quad\quad\quad\quad\quad\quad\quad\quad\quad\quad\quad\quad\quad \geq \bigg \{ \frac{b_k}{2}\int_{-\infty}^{\infty}\int_{-\infty}^{\infty}|f(\xi_1,\xi_2)|^2 d \xi_1 d \xi_2 \bigg \}^2
\end{align*}

\parindent=0mm \vspace{.1in}
Now, further applying Plancherel's theorem for the right sided BiQLCT to the right hand side of \ref{eq:4.1}, we have

\begin{align}
& \int_{-\infty}^{\infty}\int_{-\infty}^{\infty}\xi_k^2|~^{\mathbb{H}_{\mathbb{C}}}\mathcal{L}_{M_1,M_2}^{-1} [\, ^{\mathbb{H}_{\mathbb{C}}}\mathcal{L}_{M_1,M_2}^{RB}\{f\}](\xi_1,\xi_2)|^2 d \xi_1 d \xi_2 \nonumber 
 \int_{-\infty}^{\infty}\int_{-\infty}^{\infty} \omega_k^2|\,  ^{\mathbb{H}_{\mathbb{C}}}\mathcal{L}_{M_1,M_2}^{RB}\{f\}(\omega,\nu)|^2 d\omega d\nu \nonumber \\
& \quad\quad\quad\quad\quad\quad\quad\quad\quad\quad\quad\quad\quad\quad\quad\quad\quad\quad\quad \geq \bigg \{\frac{b_k}{2}\int_{-\infty}^{\infty}\int_{-\infty}^{\infty}|\,  ^{\mathbb{H}_{\mathbb{C}}}\mathcal{L}_{M_1,M_2}^{RB}\{f\}(\omega,\nu)|^2 d\omega d\nu \bigg \}^2
\end{align}

Now, we arrive at the following result.

\parindent=0mm \vspace{.1in}
{\bf{Theorem 4.2.}} (BiQWLCT uncertainity principle). For a given window function  $\phi \in L^2(\mathbb{R}^2,\mathbb{H}_{\mathbb{C}})\backslash \{0\}$, let $\mathcal{G}_{\phi}^{M_1,M_2} \{f\} \in L^2(\mathbb{R}^2,\mathbb{H}_{\mathbb{C}})$ be the BiQWLCT of $f$. Then, for every $f \in L^2(\mathbb{R}^2,\mathbb{H}_{\mathbb{C}})$, we have the following inequatlity:

\begin{align}
& \bigg\{ \int_{-\infty}^{\infty}\int_{-\infty}^{\infty} \int_{-\infty}^{\infty}\int_{-\infty}^{\infty} \omega_k^2 |\mathcal{G}_{\phi}^{M_1,M_2}f(\boldsymbol{\omega},\boldsymbol{\nu})|^2 d\omega_1 d\omega_2 d\nu_1 d\nu_2\bigg \}^{1/2} \bigg \{ \int_{-\infty}^{\infty}\int_{-\infty}^{\infty} \xi_k^2|f(\xi_1,\xi_2|^2 d\xi_1d\xi_2 \bigg\}^{1/2} \nonumber \\
& \quad\quad\quad\quad\quad\quad\quad\quad\quad\quad\quad\quad\quad\quad\quad\quad\quad\quad\quad\quad\quad \geq \frac{b_k}{4}||f||_{L^2(\mathbb{R}^2,\mathbb{H}_{\mathbb{C}})}^2||\phi||_{L^2(\mathbb{R}^2,\mathbb{H}_{\mathbb{C}})}^2
\end{align}

In order to prove theorem 4.2, we need the following lemma.

\parindent=0mm \vspace{.1in}
{\bf{Lemma 4.3.}} Suppose there exists a window function $\phi \in L^2(\mathbb{R}^2,\mathbb{H}_{\mathbb{C}})$ such that $\phi$ is not the zero function, and let $f \in L^2(\mathbb{R}^2,\mathbb{H}_{\mathbb{C}})$. Then,

\begin{align}
&||\phi||_{L^2(\mathbb{R}^2,\mathbb{H}_{\mathbb{C}})}^2\int_{-\infty}^{\infty}\int_{-\infty}^{\infty} \xi_k^2|f(\xi_1,\xi_2|^2 d\xi_1d\xi_2 \nonumber \\
&\quad\quad\quad\quad\quad\quad\quad\quad=\int_{-\infty}^{\infty}\int_{-\infty}^{\infty} \int_{-\infty}^{\infty}\int_{-\infty}^{\infty} \xi_k^2 |~ ^{\mathbb{H}_{\mathbb{C}}}\mathcal{L}_{M_1,M_2}^{-1} \{\mathcal{G}_{\phi}^{M_1,M_2}f(\boldsymbol{\omega},\boldsymbol{\nu})\}(\xi_1,\xi_2)|^2 d\xi_1 d\xi_2 d\nu_1 d\nu_2
\end{align}

{\bf{Proof.}} We have
\begin{align*}
||\phi||_{L^2(\mathbb{R}^2,\mathbb{H}_{\mathbb{C}})}^2\int_{-\infty}^{\infty}\int_{-\infty}^{\infty} \xi_k^2|f(\xi_1,\xi_2)|^2 d\xi_1d\xi_2 &= \int_{-\infty}^{\infty}\int_{-\infty}^{\infty} \xi_k^2|f(\xi_1,\xi_2|^2 d\xi_1d\xi_2 \int_{-\infty}^{\infty}\int_{-\infty}^{\infty} |\phi (\boldsymbol{\xi}-\boldsymbol{\nu})|^2d\nu_1d\nu_2 \\\\
&=\int_{-\infty}^{\infty}\int_{-\infty}^{\infty}\int_{-\infty}^{\infty}\int_{-\infty}^{\infty}\xi_k^2|f(\xi_1,\xi_2)|^2|\phi (\boldsymbol{\xi}-\boldsymbol{\nu})|^2 d\xi_1d\xi_2 d\nu_1d\nu_2\\\\
&=\int_{-\infty}^{\infty}\int_{-\infty}^{\infty}\int_{-\infty}^{\infty}\int_{-\infty}^{\infty}\xi_k^2 \big|f(\xi_1,\xi_2)\overline{\phi (\boldsymbol{\xi}-\boldsymbol{\nu})}\big|^2 d\xi_1d\xi_2 d\nu_1d\nu_2
\end{align*}

From the BiQWLCT inverse (3.3), we obtain

\begin{align}
&||\phi||_{L^2(\mathbb{R}^2,\mathbb{H}_{\mathbb{C}})}^2\int_{-\infty}^{\infty}\int_{-\infty}^{\infty} \xi_k^2|f(\xi_1,\xi_2|^2 d\xi_1d\xi_2 \nonumber \\
&\quad\quad\quad\quad\quad\quad\quad\quad=\int_{-\infty}^{\infty}\int_{-\infty}^{\infty} \int_{-\infty}^{\infty}\int_{-\infty}^{\infty} \xi_k^2 |~ ^{\mathbb{H}_{\mathbb{C}}}\mathcal{L}_{M_1,M_2}^{-1} \{\mathcal{G}_{\phi}^{M_1,M_2}f(\boldsymbol{\omega},\boldsymbol{\nu})\}(\xi_1,\xi_2)|^2 d\xi_1 d\xi_2 d\nu_1 d\nu_2
\end{align}

In the second step, we have applied Fubini's theorem to interchange the order of integration. With this, the proof is concluded.

\parindent=8mm \vspace{.1in}
Let us begin with the proof of Theorem 4.2.

\parindent=0mm \vspace{.1in}
{\bf{Proof of Theorem 4.2.}} Assume that $^{\mathbb{H}_{\mathbb{C}}}\mathcal{L}_{M_1,M_2}^{RB}\{f\} \in L^2(\mathbb{R}^2,\mathbb{H}_{\mathbb{C}}) $. Since $\mathcal{G}_{\phi}^{M_1,M_2} f \in L^2(\mathbb{R}^2,\mathbb{H}_{\mathbb{C}})$, we can replace the BiQWLCT of $f$ on both sides of (5). Then we have

\begin{align}
& \int_{-\infty}^{\infty}\int_{-\infty}^{\infty} \omega_k^2|\mathcal{G}_{\phi}^{M_1,M_2}f(\boldsymbol{\omega},\boldsymbol{\nu})|^2 d\omega_1 d\omega_2  \int_{-\infty}^{\infty}\int_{-\infty}^{\infty}\xi_k^2|~^{\mathbb{H}_{\mathbb{C}}}\mathcal{L}_{M_1,M_2}^{-1} [\, ^{\mathbb{H}_{\mathbb{C}}}\mathcal{L}_{M_1,M_2}^{RB}\{f\}](\xi_1,\xi_2)|^2 d \xi_1 d \xi_2 \nonumber \\
& \quad\quad\quad\quad\quad\quad\quad\quad\quad\quad\quad\quad\quad\quad\quad\quad\quad\quad\quad \geq \bigg \{ \frac{b_k}{2}\int_{-\infty}^{\infty}\int_{-\infty}^{\infty}|\mathcal{G}_{\phi}^{M_1,M_2}f(\boldsymbol{\omega},\boldsymbol{\nu})|^2 d\omega_1 d\omega_2 \bigg \}^2
\end{align}

Now, by applying the square root to both sides of equation (4.6) and integrating both sides with respect to $ d\boldsymbol{\nu} $, we obtain:

\begin{align}
& \int_{-\infty}^{\infty}\int_{-\infty}^{\infty}\bigg \{ \bigg (\int_{-\infty}^{\infty}\int_{-\infty}^{\infty} \omega_k^2|\mathcal{G}_{\phi}^{M_1,M_2}f(\boldsymbol{\omega},\boldsymbol{\nu})|^2 d\omega_1 d\omega_2 \bigg)^{1/2}  \nonumber \\
&   \quad\quad\quad\quad\quad\quad\quad\quad \times\bigg(\int_{-\infty}^{\infty}\int_{-\infty}^{\infty}\xi_k^2|~^{\mathbb{H}_{\mathbb{C}}}\mathcal{L}_{M_1,M_2}^{-1} [\, ^{\mathbb{H}_{\mathbb{C}}}\mathcal{L}_{M_1,M_2}^{RB}\{f\}](\xi_1,\xi_2)|^2 d \xi_1 d \xi_2 \bigg )^{1/2}\bigg\}d\nu_1 d\nu_2 \nonumber \\
& \quad\quad\quad\quad\quad\quad\quad\quad \geq  \frac{b_k}{4}\int_{-\infty}^{\infty}\int_{-\infty}^{\infty}\int_{-\infty}^{\infty}\int_{-\infty}^{\infty}|\mathcal{G}_{\phi}^{M_1,M_2}f(\boldsymbol{\omega},\boldsymbol{\nu})|^2 \, d\omega_1 \, d\omega_2 \,d\nu_1 \, d\nu_2.
\end{align}

Additionally, by applying the Cauchy-Schwarz inequality to the left-hand side of equation (4.7), we can deduce:

\begin{align}
&\bigg\{\int_{-\infty}^{\infty}\int_{-\infty}^{\infty} \int_{-\infty}^{\infty}\int_{-\infty}^{\infty} \omega_k^2|\mathcal{G}_{\phi}^{M_1,M_2}f(\boldsymbol{\omega},\boldsymbol{\nu})|^2 \, d\omega_1 \, d\omega_2 \,d\nu_1 \, d\nu_2\bigg\}^{1/2}  \nonumber \\\
&\quad\quad\quad\quad\quad\quad \times \bigg\{\int_{-\infty}^{\infty}\int_{-\infty}^{\infty} \int_{-\infty}^{\infty}\int_{-\infty}^{\infty}\xi_k^2|~^{\mathbb{H}_{\mathbb{C}}}\mathcal{L}_{M_1,M_2}^{-1} [\, ^{\mathbb{H}_{\mathbb{C}}}\mathcal{L}_{M_1,M_2}^{RB}\{f\}](\xi_1,\xi_2)|^2 \, d\xi_1 \, d\xi_2 \, d\nu_1 \, d\nu_2 \bigg\}^{1/2} \nonumber \\\nonumber\\
&\qquad\qquad\quad\qquad\qquad\quad\quad\quad \geq  \frac{b_k}{4}\int_{-\infty}^{\infty}\int_{-\infty}^{\infty}\int_{-\infty}^{\infty}\int_{-\infty}^{\infty}|\mathcal{G}_{\phi}^{M_1,M_2}f(\boldsymbol{\omega},\boldsymbol{\nu})|^2 \, d\omega_1 \, d\omega_2 \, d\nu_1 \, d\nu_2
\end{align}

Inserting Lemma 4.2 into the second term on the left-hand side of (4.8) and substituting 
$$\int_{-\infty}^{\infty}\int_{-\infty}^{\infty}\int_{-\infty}^{\infty}\int_{-\infty}^{\infty}\big | \mathcal{G}_{\phi}^{M_1,M_2}f(\boldsymbol{\omega},\boldsymbol{\nu})\big|^2 \, d\omega_1 \, d\omega_2 \, d\nu_1 \, d\nu_2 =||f||_{L^2(\mathbb{R}^2,\mathbb{H}_{\mathbb{C}})}^2||\phi||_{L^2(\mathbb{R}^2,\mathbb{H}_{\mathbb{C}})}^2$$
 into the right-hand side of this inequality. Then, we have

\begin{align}
& \bigg\{\int_{-\infty}^{\infty}\int_{-\infty}^{\infty} \int_{-\infty}^{\infty}\int_{-\infty}^{\infty} \omega_k^2|\mathcal{G}_{\phi}^{M_1,M_2}f(\boldsymbol{\omega},\boldsymbol{\nu})|^2 \, \, d\omega_1 \, d\omega_2 \,d\nu_1 \, d\nu_2\bigg\}^{1/2} \nonumber \\
&\quad\quad\quad\quad
\quad \times \bigg\{||\phi||_{L^2(\mathbb{R}^2,\mathbb{H}_{\mathbb{C}})} \int_{-\infty}^{\infty}\int_{-\infty}^{\infty} \int_{-\infty}^{\infty}\int_{-\infty}^{\infty} \xi_k^2|f(\xi_1,\xi_2|^2 d\xi_1d\xi_2 \bigg\}^{1/2} \nonumber \\\
&\quad\quad\quad\quad
\quad \ge  \frac{b_k}{4}||f||_{L^2(\mathbb{R}^2,\mathbb{H}_{\mathbb{C}})}^2||\phi||_{L^2(\mathbb{R}^2,\mathbb{H}_{\mathbb{C}})}^2.
\end{align}
Dividing both sides of (4.9) by $||\phi||_{L^2(\mathbb{R}^2,\mathbb{H}_{\mathbb{C}})}$, we obtain

\begin{align}
& \bigg\{\int_{-\infty}^{\infty}\int_{-\infty}^{\infty} \int_{-\infty}^{\infty}\int_{-\infty}^{\infty} \omega_k^2 |\mathcal{G}_{\phi}^{M_1,M_2}f(\boldsymbol{\omega},\boldsymbol{\nu})|^2 \, d\omega_1 \, d\omega_2 \,d\nu_1 \, d\nu_2\bigg \}^{1/2} \bigg \{ \int_{-\infty}^{\infty}\int_{-\infty}^{\infty} \xi_k^2|f(\xi_1,\xi_2|^2 d\xi_1d\xi_2 \bigg\}^{1/2} \nonumber \\\\
& \quad\quad\quad\quad\quad\quad\quad\quad\quad\quad\quad\quad\quad\quad\quad\quad\quad\quad\quad\quad\quad \geq \frac{b_k}{4}||f||_{L^2(\mathbb{R}^2,\mathbb{H}_{\mathbb{C}})}^2||\phi||_{L^2(\mathbb{R}^2,\mathbb{H}_{\mathbb{C}})}
\end{align}
which completes the proof.

\parindent=0mm \vspace{.3in}
 { \bf{5. Example of the right sided BiQWLCT}} 

\parindent=0mm \vspace{.1in}
Considering the window function associated with the two-dimensional Haar function as described by
$$\phi(\boldsymbol{\xi})=\left\{ \begin{array}{ll}
1, & \text{if } 0 \leq \xi_1 < \frac{1}{2}, \, 0 \leq \xi_2 < \frac{1}{2}\\
-1, & \text{if } \frac{1}{2} \leq \xi_1 < 1, \, \frac{1}{2} \leq \xi_2 < 1\\
0, & \text{Otherwise}
\end{array}\right.$$

Examine a  biquaternionic Gaussian function $f$ represented by the following expression 
$$f(\xi_1,\xi_2)=C_0e^{-(\alpha_1\xi_1^2+\alpha_2\xi_2^2)}$$
where $C_0$ is a complex constant and $\alpha_1, \alpha_2 \in \mathbb{R}$ are positive real constants.

\parindent=0mm \vspace{.1in}
Then the right sided BiQWLCT of $f$ is given by
\begin{align*}
\mathcal{G}_{\phi}^{M_1,M_2}\{f(\boldsymbol{\xi})\}( \boldsymbol{\omega},\boldsymbol{\nu})&=\int_{-\infty}^{\infty}\int_{-\infty}^{\infty}f(\boldsymbol{\xi})\overline{\phi(\boldsymbol{\xi}-\boldsymbol{\nu})}K_{A_1}^{\boldsymbol{\mu}}(\xi_1,\omega_1) K_{A_2}^{\boldsymbol{\mu}}(\xi_2,\omega_2)d\xi_1d\xi_2\\\\
&=\int_{-\infty}^{\infty}\int_{-\infty}^{\infty}f(\boldsymbol{\xi})\overline{\phi(\boldsymbol{\xi}-\boldsymbol{\nu})}\frac{1}{\sqrt{2\pi|b_1|}} e^{\boldsymbol{\mu} \big (\frac{a_1}{2b_1}\xi_1^2-\frac{\xi_1\omega_1}{b_1}+\frac{d_1}{2b_1}\omega_1^2-\frac{\pi}{4}\big )}\\
& \qquad\qquad\qquad\qquad\times \frac{1}{\sqrt{2\pi|b_2|}} e^{\boldsymbol{\mu} \big (\frac{a_2}{2b_2}\xi_2^2-\frac{\xi_2\omega_2}{b_2}+\frac{d_2}{2b_2}\omega_2^2-\frac{\pi}{4}\big )} d\xi_1 d\xi_2\\
&= \frac{C_0}{2\pi\sqrt{|b_1b_2|}} \int_{\omega_1}^{1/2+\omega_1}e^{-\alpha_1\xi_1^2}e^{\boldsymbol{\mu} \big (\frac{a_1}{2b_1}\xi_1^2-\frac{\xi_1\omega_1}{b_1}+\frac{d_1}{2b_1}\omega_1^2-\frac{\pi}{4}\big )}d\xi_1\\
& \qquad\qquad\qquad\qquad\times \int_{\omega_2}^{1/2+\omega_2}e^{-\alpha_2\xi_2^2}e^{\boldsymbol{\mu} \big (\frac{a_2}{2b_2}\xi_2^2-\frac{\xi_2\omega_2}{b_2}+\frac{d_2}{2b_2}\omega_2^2-\frac{\pi}{4}\big )}  d\xi_2\\\\
&\qquad- \frac{C_0}{2\pi\sqrt{|b_1b_2|}} \int_{1/2+\omega_1}^{1+\omega_1}e^{-\alpha_1\xi_1^2} e^{\boldsymbol{\mu} \big (\frac{a_1}{2b_1}\xi_1^2-\frac{\xi_1\omega_1}{b_1}+\frac{d_1}{2b_1}\omega_1^2-\frac{\pi}{4}\big )}d\xi_1\\\\
&\qquad\qquad\qquad\qquad \times \int_{1/2+\omega_1}^{1+\omega_1}e^{-\alpha_2\xi_2^2} e^{\boldsymbol{\mu} \big (\frac{a_2}{2b_2}\xi_2^2-\frac{\xi_2\omega_2}{b_2}+\frac{d_2}{2b_2}\omega_2^2-\frac{\pi}{4}\big )}d\xi_2\\
\end{align*}

By using the method of completing squares, we obtain
\begin{align*}
&\mathcal{G}_{\phi}^{M_1,M_2}\{f(\boldsymbol{\xi})\}( \boldsymbol{\omega},\boldsymbol{\nu})\\
&\quad \quad \quad \quad= \frac{C_0}{2\pi\sqrt{|b_1b_2|}} \int_{\omega_1}^{1/2+\omega_1} e^{\big( \frac{-\omega_1^2}{2b_1\alpha_1\xi_1-\boldsymbol{\mu}\alpha_1} +\boldsymbol{\mu}\frac{d_1}{2b_1}\omega_1^2-\frac{\pi}{4} \big)}\cdot e^{-\big ( \sqrt{\frac{2b_1\alpha_1\xi_1-\boldsymbol{\mu}\alpha_1}{2b_1}}\xi_1+\frac{\boldsymbol{\mu}\omega_1}{\sqrt{2b_1(2b_1\alpha_1\xi_1-\boldsymbol{\mu}\alpha_1)}} \big )^2}d\xi_1\\\\
&\qquad\qquad\quad\quad \quad \quad \quad \times   \int_{\omega_2}^{1/2+\omega_2} e^{\big( \frac{-\omega_2^2}{2b_2\alpha_2\xi_2-\boldsymbol{\mu}\alpha_2} +\boldsymbol{\mu}\frac{d_2}{2b_2}\omega_2^2-\frac{\pi}{4} \big)}\cdot e^{-\big ( \sqrt{\frac{2b_2\alpha_2\xi_2-\boldsymbol{\mu}\alpha_2}{2b_2}}\xi_2+\frac{\boldsymbol{\mu}\omega_2}{\sqrt{2b_2(2b_2\alpha_2\xi_2-\boldsymbol{\mu}\alpha_2)}} \big )^2}d\xi_2\\\\
&\quad \quad\qquad \quad \quad- \frac{C_0}{2\pi\sqrt{|b_1b_2|}} \int_{1/2+\omega_1}^{1+\omega_1} e^{\big( \frac{-\omega_1^2}{2b_1\alpha_1\xi_1-\boldsymbol{\mu}\alpha_1} +\boldsymbol{\mu}\frac{d_1}{2b_1}\omega_1^2-\frac{\pi}{4} \big)}\cdot e^{-\big ( \sqrt{\frac{2b_1\alpha_1\xi_1-\boldsymbol{\mu}\alpha_1}{2b_1}}\xi_1+\frac{\boldsymbol{\mu}\omega_1}{\sqrt{2b_1(2b_1\alpha_1\xi_1-\boldsymbol{\mu}\alpha_1)}} \big )^2}d\xi_1\\\\
&\quad \qquad\qquad\qquad\quad \quad \quad \times   \int_{1/2+\omega_2}^{1+\omega_2} e^{\big( \frac{-\omega_2^2}{2b_2\alpha_2\xi_2-\boldsymbol{\mu}\alpha_2} +\boldsymbol{\mu}\frac{d_2}{2b_2}\omega_2^2-\frac{\pi}{4} \big)}\cdot e^{-\big ( \sqrt{\frac{2b_2\alpha_2\xi_2-\boldsymbol{\mu}\alpha_2}{2b_2}}\xi_2+\frac{\boldsymbol{\mu}\omega_2}{\sqrt{2b_2(2b_2\alpha_2\xi_2-\boldsymbol{\mu}\alpha_2)}} \big )^2}d\xi_2
\end{align*}

Now making the substitutions $M_1=\sqrt{\frac{2b_1\alpha_1\xi_1-\boldsymbol{\mu}\alpha_1}{2b_1}}$, $M_2=\sqrt{\frac{2b_2\alpha_2\xi_2-\boldsymbol{\mu}\alpha_2}{2b_2}}$,

\parindent=0mm \vspace{.in}
$N_1=\frac{\boldsymbol{\mu}\omega_1}{\sqrt{2b_1(2b_1\alpha_1\xi_1-\boldsymbol{\mu}\alpha_1)}}$,  $N_2=\frac{\boldsymbol{\mu}\omega_2}{\sqrt{2b_2(2b_2\alpha_2\xi_2-\boldsymbol{\mu}\alpha_2)}}$,  $Q_1=e^{\big( \frac{-\omega_1^2}{2b_1\alpha_1\xi_1-\boldsymbol{\mu}\alpha_1} +\boldsymbol{\mu}\frac{d_1}{2b_1}\omega_1^2-\frac{\pi}{4} \big)}$ and $Q_2=e^{\big( \frac{-\omega_2^2}{2b_2\alpha_2\xi_2-\boldsymbol{\mu}\alpha_2} +\boldsymbol{\mu}\frac{d_2}{2b_2}\omega_2^2-\frac{\pi}{4} \big)}$ in the above expression we obtain

\begin{align*}
\mathcal{G}_{\phi}^{M_1,M_2}\{f(\boldsymbol{\xi})\}( \boldsymbol{\omega},\boldsymbol{\nu})&= \frac{C_0}{2\pi\sqrt{|b_1b_2|}} \int_{\omega_1}^{1/2+\omega_1} Q_1 e^{-(M_1\xi_1+N_1)^2}d\xi_1 \times \int_{\omega_2}^{1/2+\omega_2} Q_2 e^{-(M_2\xi_2+N_2)^2}d\xi_2\\\\
&\qquad- \frac{C_0}{2\pi\sqrt{|b_1b_2|}} \int_{1/2+\omega_1}^{1+\omega_1} Q_1 e^{-(M_1\xi_1+N_1)^2}d\xi_1 \times \int_{1/2+\omega_2}^{1+\omega_2} Q_2 e^{-(M_2\xi_2+N_2)^2}d\xi_2
\end{align*}

Now, substituting $\nu_1=M_1\xi_1+N_1$ and $\nu_2=M_2\xi_2+N_2$ in the above equation, we have

\begin{align*}
\mathcal{G}_{\phi}^{M_1,M_2}\{f(\boldsymbol{\xi})\}( \boldsymbol{\omega},\boldsymbol{\nu})&= \frac{C_0Q_1	Q_2}{M_1M_2 2\pi\sqrt{|b_1b_2|}} \int_{M_1\omega_1+N_1}^{M_1(1/2+\omega_1)+N_1}e^{-\nu_1^2}d\nu_1 \times \int_{M_2\omega_2+N_2}^{M_2(1/2+\omega_2)+N_2}e^{-\nu_2^2}d\nu_2 \\\\
&\qquad- \frac{C_0Q_1	Q_2}{M_1M_2 2\pi\sqrt{|b_1b_2|}} \int_{M_1(1/2+\omega_1)+N_1}^{M_1(1+\omega_1)+N_1}e^{-\nu_1^2}d\nu_1 \times \int_{M_2(1/2+\omega_2)+N_2}^{M_2(1+\omega_2)+N_2}e^{-\nu_2^2}d\nu_2
\end{align*}
This can be further written as

\begin{align*}
\mathcal{G}_{\phi}^{M_1,M_2}\{f(\boldsymbol{\xi})\}( \boldsymbol{\omega},\boldsymbol{\nu})&= \frac{C_0Q_1	Q_2}{M_1M_2 2\pi\sqrt{|b_1b_2|}} \bigg\{\int_{0}^{M_1\omega_1+N_1}-e^{-\nu_1^2}d\nu_1 + \int_{0}^{M_1(1/2+\omega_1)+N_1}e^{-\nu_1^2}d\nu_1 \bigg\} \\\\
&\qquad\qquad\qquad\quad \times \bigg\{\int_{0}^{M_2\omega_2+N_2}(-e^{-\nu_2^2})d\nu_2 + \int_{0}^{M_2(1/2+\omega_2)+N_2}e^{-\nu_2^2}d\nu_2 \bigg\}\\\\
&\quad- \frac{C_0Q_1	Q_2}{M_1M_2 2\pi\sqrt{|b_1b_2|}} \bigg\{\int_{0}^{M_1(1/2+\omega_1)+N_1}-e^{-\nu_1^2}d\nu_1 + \int_{0}^{M_1(1+\omega_1)+N_1}e^{-\nu_1^2}d\nu_1 \bigg\} \\\\
&\qquad\qquad\qquad \times \bigg\{\int_{0}^{M_2(1/2+\omega_2)+N_2}(-e^{-\nu_2^2})d\nu_2 + \int_{0}^{M_2(1+\omega_2)+N_2}e^{-\nu_2^2}d\nu_2 \bigg\}. \square
\end{align*}

\parindent=0mm \vspace{.3in}
 { \bf{6. Potential Applications}}

\parindent=0mm \vspace{.1in}
Due to its capability of analyzing complex signals and images in both time and frequency domains at the same time, Biquaternion Windowed Linear Canonical Transform (BiQWLCT) has a number of uses in signal processing as well as image processing. In this regard, some of its possible applications include:

\parindent=0mm \vspace{.1in}

{\it 1.Signal Analysis}: 

\parindent=0mm \vspace{.1in}
{\it Time-Frequency Analysis:} BiQWLCT allows simultaneous analysis of signals in the time and frequency domains, providing a detailed representation of signal characteristics.BiQWLCT offers advantages over traditional methods like Fourier transforms by providing time-frequency analysis along with feature extraction. This allows for a more detailed understanding of the signal's behavior in different time intervals and frequencies.

\parindent=0mm \vspace{.1in}
{\it Modulation and Demodulation:}  BiQWLCT can be used for efficient modulation and demodulation of communication signals in systems where signal components are distributed over a wide frequency range.

\parindent=0mm \vspace{.1in}
{\it Feature Extraction: }In pattern recognition tasks, BiQWLCT can extract discriminative features from signals for classification and identification purposes.

\parindent=0mm \vspace{.1in}
{\it 2. Image Processing:}

\parindent=0mm \vspace{.1in}
{\it Image Enhancement:} BiQWLCT can enhance images by analyzing and manipulating their frequency components, leading to improved visual quality and detail.

\parindent=0mm \vspace{.1in}
{\it Texture Analysis:} By extracting texture features using BiQWLCT, images can be analyzed for patterns, textures, and structures that are useful in various image processing tasks.

\parindent=0mm \vspace{.1in}
{\it Image Compression:} BiQWLCT can be utilized for image compression by analyzing and representing image data efficiently in the transformed domain.

\parindent=8mm \vspace{.1in}
Overall, BiQWLCT is a promising tool for advanced signal and image processing tasks, particularly when dealing with multidimensional data and extracting local features.

\parindent=0mm\vspace{0.2in}
{\bf{Acknowledgment}}

\parindent=0mm\vspace{0.in}
Second author acknowledges Department of Science and Technology (DST) ,Govt. of India for financial support through INSPIRE Fellowship(grant no:DST/INSPIRE/2021/IF210207).

\parindent=0mm\vspace{0.2in}
{\bf{Data Availability}} Data sharing not applicable to this article as no datasets were generated or analysed during the current study.

\parindent=0mm\vspace{0.2in}
{\bf{Declarations}}

\parindent=0mm\vspace{0.in}
{\bf{Conflict of Interest}} The author declares that he has no conflict of interest.

\parindent=0mm\vspace{0.5in}
{\bf{\LARGE{References}}}

\end{document}